\documentstyle[12pt]{article}
\newtheorem{thm}{Theorem}

\newtheorem{remark}{Remark}
\newtheorem{lemma}{Lemma}

\newtheorem{definition}{Definition}
\newtheorem{cor}{Corollary}

\newcommand{\pfbegin}{\noindent {\bf Proof:\ }}
\newcommand{\pfend}{\hspace*{\fill}{\bf $\Box$} \bigskip \\}

\newcommand{\N}{I\!\!N}

\newcommand{\de}{\delta}
\newcommand{\g}{\gamma}
\title{An Algorithmic Proof of Suslin's Stability Theorem over
Polynomial Rings}
\author{Hyungju Park\thanks{Dept. of Mathematics, University of
California, Berkeley; park@math.berkeley.edu} \and
Cynthia Woodburn\thanks{Dept. of Mathematics, Pittsburg State
University; cwoodbur@mail.pittstate.edu}}
\date{}
\begin{document}
\maketitle
\begin{abstract}
Let $k$ be a field. Then Gaussian elimination over $k$ and the Euclidean
division algorithm for the univariate polynomial ring $k[x]$
allow us to write any matrix in
$SL_n(k)$ or $SL_n(k[x])$, $n\geq 2$, as a product of elementary
matrices. Suslin's stability theorem states that the same is true for
the multivariate polynomial ring $SL_n(k[x_1,\ldots ,x_m])$
with $n\geq 3$. As
Gaussian elimination gives us an algorithmic way of finding an explicit
factorization of the given matrix into elementary matrices over a field,
we develop a similar algorithm over polynomial rings.
\end{abstract}
%
%This is Section 1
%
\section{Introduction}
Immediately after proving the famous {\em Serre's Conjecture}
(the {\em Quillen-Suslin theorem}, nowadays) in 1976 \cite{suslin:serre},
A. Suslin went on \cite{suslin} to prove the following $K_1$-analogue of
{\em Serre's Conjecture} which is now known as
{\em Suslin's stability theorem}:
\begin{quote}
Let $R$ be a commutative Noetherian ring and $n\geq \max(3,\dim(R)+2)$.
Then, any $n\times n$ matrix $A=(f_{ij})$ of determinant $1$, with
$f_{ij}$ being elements of the polynomial ring $R[x_1,\ldots ,x_m]$, can be
written as a product of elementary matrices over $R[x_1,\ldots ,x_m]$.
\end{quote}
\begin{definition}
For any ring $R$, an $n\times n$ elementary matrix $E_{ij}(a)$ over $R$ is
a matrix of the form
$I+a\cdot e_{ij}$ where $i\neq j,a\in R$ and $e_{ij}$ is the
$n\times n$ matrix whose $(i,j)$ component is $1$ and all other
components are zero.
\end{definition}
For a ring $R$, let $SL_n(R)$ be be the group of all the $n\times n$
matrices of determinant $1$ whose entries are elements of $R$, and let
$E_n(R)$ be the subgroup of $SL_n(R)$ generated by the elementary matrices.
Then {\em Suslin's stability theorem} can be expressed as
\begin{eqnarray}
	SL_n(R[x_1,\ldots ,x_m]) =  E_n(R[x_1,\ldots ,x_m])
\quad {\rm for\ all}\ n\geq \max(3,\dim(R)+2).
\end{eqnarray}
In this paper, we develop an algorithmic proof of the above assertion
over a field $k$. By implementing this algorithm, for a given
$A\in SL_n(k[x_1,\ldots ,x_m])$
with $n\geq 3$, we are able to find those elementary
matrices $E_1,\ldots ,E_t\in E_n(k[x_1,\ldots ,x_m])$ such that
$A=E_1\cdots E_t$.
\begin{remark}
If a matrix $A$ can be written as a product of elementary matrices, we
will say $A$ is {\em realizable}.
\end{remark}
\bigskip
\begin{itemize}
\item
In section 2, an algorithmic proof of the normality of
$E_n(k[x_1,\ldots ,x_m])$ in $SL_n(k[x_1,\ldots ,x_m])$ for $n\geq 3$
is given, which will be used in the rest of paper.
\item
In section 3, we develop an algorithm for the {\em Quillen Induction Process},
a standard way of reducing a given problem over a ring to an
easier problem over a local ring. Using this {\em Quillen Induction Algorithm},
we reduce our realization problem over the polynomial
ring $R[X]$ to one over $R_M[X]$'s, where $R=k[x_1,\ldots, x_{m-1}]$
and $M$ is a maximal ideal of $R$.
\item
In section 4, an algorithmic proof of the {\em Elementary Column Property},
a stronger version of the {\em Unimodular Column Property}, is given,
and we note that this algorithm gives  another constructive proof
of the {\em Quillen-Suslin theorem}.
Using the {\em Elementary Column Property}, we show that a
realization algorithm for
$SL_n(k[x_1,\ldots ,x_m])$ is obtained from
a realization algorithm for the matrices of the following
special form:
$$\left(\! \begin{array}{ccc} p & q & 0 \\ r & s & 0 \\ 0 & 0 & 1
\end{array}\!\right) \in SL_3(k[x_1,\ldots ,x_m]),$$
where $p$ is monic in the last variable $x_m$.
\item
In section 5, in view of the results in the preceding two sections, we note
that a realization algorithm over $k[x_1,\ldots ,x_m]$ can be obtained
from a realization algorithm for the matrices of the special
form $\left(\!\begin{array}{ccc} p & q & 0 \\ r & s & 0 \\ 0 & 0 & 1
\end{array}\!\right)$
over $R[X]$, where $R$ is now a local ring and $p$ is monic in $X$.
A realization algorithm
for this case was already found by M.P. Murthy in \cite{murthy}. We
reproduce {\em Murthy's Algorithm} in this section.
\item
In section 6, we suggest using the {\em Steinberg relations} from algebraic
$K$-theory to lower the number of elementary matrix factors in a
factorization produced by our algorithm.
We also mention an ongoing effort of using our algorithm in {\em Signal
Processing}.

\end{itemize}
%
%
%This is Section 2
%
\section{Normality of $E_n(k[x_1,\ldots ,x_m])$ in $SL_n(k[x_1,\ldots ,x_m])$}
\begin{lemma}
The Cohn matrix $A=\left(\! \begin{array}{cc} 1+xy & x^2 \\ -y^2 & 1-xy
\end{array}\!\right)$ is not realizable, but $\left(\! \begin{array}{cc}
A & 0 \\ 0 & 1 \end{array}\!\right) \in SL_3(k[x,y])$ is.
\end{lemma}
\pfbegin
The nonrealizability of $A$ is proved in \cite{cohn}, and a complete
algorithmic criterion for the realizability of
matrices in $SL_2(k[x_1,\ldots ,x_m])$
is developed in \cite{thk}. Now consider
\begin{eqnarray}
\left(\! \begin{array}{cc} A & 0 \\ 0 & 1 \end{array}
\!\right) =\left(\! \begin{array}{ccc} 1+xy & x^2 & 0 \\ -y^2 & 1-xy & 0
\\ 0 & 0 & 1 \end{array}\!\right).
\end{eqnarray}
Noting that $\left(\!\begin{array}{ccc} 1+xy & x^2 & 0 \\ -y^2 & 1-xy & 0
\\ 0 & 0 & 1 \end{array}\!\right)=I+\left(\! \begin{array}{c} x\\ -y \\ 0
\end{array}\!\right)\cdot (y,x,0)$, we see that
the realizability of this matrix is
a special case of the following {\bf Lemma 3}.
\pfend
\begin{definition}
Let $n\geq 2$. A {\bf Cohn-type matrix} is a matrix of the form
$$ I+a{\bf v} \cdot (v_j{\bf e_i}-v_i{\bf e_j}) $$
where ${\bf v}=\left(\!\begin{array}{c}v_1\\ \vdots\\ v_n\end{array}\!
\right)\in (k[x_1,\ldots ,x_m])^n$, $i<j\in\{ 1,\ldots ,n\}$,
$a\in k[x_1,\ldots ,x_m]$, and ${\bf e_i}=(0,\ldots ,0,1,0,\ldots,0)$
with $1$ occurring only at the $i$-th position.
\end{definition}
\begin{lemma} \label{lem;Cohn}
Any Cohn-type matrix for $n\geq 3$ is realizable.
\end{lemma}
\pfbegin
First, let's consider the case $i=1,j=2$.
In this case,
\begin{eqnarray}
B & = & I+a\left(\!\begin{array}{c} v_1\\ \vdots\\ v_n \end{array}\!\right)
		\cdot (v_2, -v_1,0,\ldots , 0) \nonumber\\
	& = & \left(\!\begin{array}{ccccc}
		1+av_1v_2 & -av_1^2 & 0 & \cdots & 0 \\
		av_2^2 & 1-av_1v_2 & 0 & \cdots & 0 \\
		av_3v_2 & -av_3v_1 & & & \\
		\vdots & \vdots & & I_{n-2} & \\
		av_nv_2 & -av_nv_1 & & & \end{array}\!\right) \nonumber\\
	& = & \left(\!\begin{array}{ccccc}
		1+av_1v_2 & -av_1^2 & 0 & \cdots & 0 \\
		av_2^2 & 1-av_1v_2 & 0 & \cdots & 0 \\
		0 & 0 & & & \\
		\vdots & \vdots & & I_{n-2} & \\
		0 & 0 & & & \end{array}\!\right)
		\prod_{l=3}^nE_{l1}(av_lv_2)E_{l2}(-av_lv_1),
\end{eqnarray}
So, it's enough to show that
\begin{eqnarray}
A=\left( \!\begin{array}{ccc} 1+av_1v_2 & -av_1^2 & 0 \\
		av_2^2 & 1-av_1v_2 & 0 \\
		0 & 0 & 1 \end{array}\!\right)
\end{eqnarray}
is realizable for any $a,v_1,v_2\in  k[x_1,\ldots ,x_m]$.
Let ``$\rightarrow$'' indicate that we are applying elementary operations,
and consider the following:
\begin{eqnarray}
A & = & \left(\! \begin{array}{ccc} 1+av_1v_2 & -av_1^2 & 0 \\
	av_2^2 & 1-av_1v_2 & 0 \\ 0 & 0 & 1\end{array}\!\right)
	\rightarrow
	\left(\!\begin{array}{ccc} 1+av_1v_2 & -av_1^2 & v_1 \\
	av_2^2 & 1-av_1v_2 & v_2\\ 0 & 0 & 1\end{array}\!\right)\nonumber\\
	& \rightarrow &
	\left( \!\begin{array}{ccc} 1 & -av_1^2 & v_1\\
	0 & 1-av_1v_2 & v_2 \\ -av_2 & 0 & 1\end{array}\!\right)
        \rightarrow
	\left(\! \begin{array}{ccc} 1 & 0 & v_1 \\ 0 & 1 & v_2 \\
	-av_2 & av_1 & 1 \end{array} \!\right)
	\rightarrow
	\left( \!\begin{array}{ccc} 1 & 0 & v_1 \\ 0 & 1 & v_2 \\
	0 & av_1 & 1+av_1v_2 \end{array}\!\right)\nonumber\\
	& \rightarrow &
	\left(\! \begin{array}{ccc} 1 & 0 & 0 \\ 0 & 1 & v_2 \\
	0 & av_1 & 1+av_1v_2 \end{array}\!\right)
        \rightarrow
	\left(\! \begin{array}{ccc} 1 & 0 &0 \\ 0 & 1 & v_2 \\
	0 & 0 & 1 \end{array}\! \right)
	\rightarrow
	\left(\!\begin{array}{ccc}
	1 & 0 &0 \\ 0 & 1 &0 \\ 0 & 0 & 1\end{array}\! \right).
\end{eqnarray}
Keeping track of all the elementary operations involved, we get
\begin{eqnarray}
A=E_{13}(-v_1)E_{23}(-v_2)E_{31}(-av_2)E_{32}(av_1)E_{13}(v_1)
E_{23}(v_2)E_{31}(av_2)E_{32}(-av_1).
\end{eqnarray}
In general (i.e., for arbitrary $i<j$),
\begin{eqnarray}
B & = & I+a\left(\!\begin{array}{c} v_1 \\ \vdots \\v_n \end{array}\!\right)
	\cdot (0,\ldots ,0,v_j,0,\ldots ,0,-v_i,0,\ldots ,0) \nonumber\\
	& & \mbox{(Here,\ $v_j$\ occurs\ at\ the\ $i$-th\ position\
	and\ $-v_i$\ occurs\ at\ the}\nonumber\\
	& & j\mbox{-th\ position.)}\nonumber\\
	& = & \left(\!\begin{array}{ccccccc}
		1 & \cdots & av_1v_j & \cdots & -av_1v_i & \cdots & 0 \\
		& \ddots & \vdots & & \vdots & & 0 \\
		& & 1+av_iv_j & & -av_i^2 & & \\
		& & \vdots & & \vdots & & \\
		& & av_j^2 & & 1-av_iv_j & & \\
		& & \vdots & & \vdots & & \\
		& & v_nv_j & & -v_nv_i & & 1 \end{array}\!\right)\nonumber\\
	& = & \left(\!\begin{array}{ccccccc}
		1 & \cdots & 0 & \cdots & 0 & \cdots & 0 \\
		& \ddots & \vdots & & \vdots & & 0 \\
		& & 1+av_iv_j & & -av_i^2 & & \\
		& & \vdots & & \vdots & & \\
		& & av_j^2 & & 1-av_iv_j & & \\
		& & \vdots & & \vdots & & \\
		& & 0 & & 0 & & 1 \end{array}\!\right)\nonumber\\
	& & \cdot\prod_{1\leq l\leq n, l\neq i,j}
		E_{li}(av_lv_j)E_{lj}(-av_lv_i)
		\nonumber\\
	& = & E_{it}(-v_i)E_{jt}(-v_j)E_{ti}(-av_j)E_{tj}(av_i)E_{it}(v_i)
		E_{jt}(v_j)E_{ti}(av_j)E_{tj}(-av_i) \nonumber\\
	& & \cdot\prod_{1\leq l\leq n, l\neq i,j}
		E_{li}(av_lv_j)E_{lj}(-av_lv_i).
\end{eqnarray}
In the above, $t\in \{1,\ldots, n\}$ can be chosen to be any number other
than $i$ and $j$.
\pfend
Since a {\em Cohn-type matrix} is realizable, any product of
{\em Cohn-type matrices} is also realizable.
This observation motivates the following generalization of the above lemma.
\begin{definition}
Let $R$ be a ring and ${\bf v}=(v_1,\ldots ,v_n)^t
\in R^n$ for some $n\in \N$. Then ${\bf v}$ is called a {\em unimodular
column vector} if its components generate $R$, i.e.
if there exist $g_1,\ldots ,g_n\in R$ such that
$v_1g_1+\cdots +v_ng_n=1$.
\end{definition}
\begin{cor}
Suppose that $A \in SL_n(k[x_1, \ldots ,x_m])$ with $n\geq 3$
can be written in the
form $A=I+{\bf v}\cdot {\bf w}$ for a unimodular column vector $\ {\bf v}$
and a row vector $\ {\bf w}$ over $\ k[x_1,\ldots ,x_m]\ $
such that ${\bf w}\cdot {\bf v}=0$. Then $A$ is realizable.
\end{cor}
\pfbegin
Since ${\bf v}=(v_1,\ldots ,v_n)^t$ is unimodular, we
can find $g_1,\ldots ,g_n \in k[x_1,\ldots ,x_m]$ such that
$v_1g_1+\cdots +v_ng_n=1$. We can use the {\em effective Nullstellensatz}
to explicitly find these $g_i$'s (See \cite{fitchas:galligo}).
This combined with ${\bf w}\cdot {\bf v}=w_1v_1+\cdots +w_nv_n=0$
yields a new expression for ${\bf w}$:
\begin{eqnarray}
     {\bf w}=\sum_{i<j}a_{ij}(v_j{\bf e_i}-v_i{\bf e_j})
\end{eqnarray}
where $a_{ij}=w_ig_j-w_jg_i$. Now,
\begin{eqnarray}
     A=\prod_{i<j}\left(I+{\bf v}\cdot a_{ij}(v_j{\bf e_i}-v_i{\bf e_j})
     \right).
\end{eqnarray}
Each component on the right hand side of this equation is a {\em Cohn-type
matrix} and thus realizable, so $A$ is also realizable.
\pfend
\begin{cor}
\quad $BE_{ij}(a)B^{-1}$  is realizable for any $\ B\in GL_n(k[x_1, \ldots
,x_m])$ with $n\geq 3$ and $\ a\in k[x_1, \ldots ,x_m]$.
\end{cor}
\pfbegin
Note that $i\neq j$, and
$$BE_{ij}(a)B^{-1}=I+(i\mbox{-th\ column\ vector\ of}\ B)\cdot a
\cdot (j\mbox{-th\ row\ vector\ of}\ B^{-1}).$$
Let ${\bf v}$ be the $i$-th column vector of $B$ and ${\bf w}$ be
$a$ times the $j$-th row vector of $B^{-1}$. Then $(i$-th row vector of
$B^{-1})\cdot {\bf
v}=1$ implies ${\bf v}$ is unimodular, and ${\bf w}\cdot {\bf v}$ is
clearly zero since $i\neq j$. Therefore, $BE_{ij}(a)B^{-1}=I+{\bf v}\cdot
{\bf w}$ satisfies the condition of the above corollary, and is thus
realizable.
\pfend
\begin{remark}
One important consequence of this corollary is that for $n\geq 3$,
$E_n(k[x_1,\ldots ,x_m])$ is a normal subgroup of $SL_n(k[x_1,\ldots ,x_m])$,
i.e. if $A\in SL_n(k[x_1,\ldots ,x_m])$ and $E\in E_n(k[x_1,\ldots ,x_m])$,
then the above corollary gives us an algorithm for finding elementary matrices
$E_1,\ldots ,E_t$ such that $A^{-1}EA=E_1\cdots E_t$.
\end{remark}
%%
%This is about Glueing Algorithm
%
\section{Glueing of Local Realizability}
Let $R=k[x_1,\ldots , x_{m-1}],\  X=x_m$ and $M\in$ Max($R$)
=\{ maximal ideals of $R$\}.
For $A\in SL_n(R[X])$, we let $A_M\in SL_n(R_M[X])$ be
its image under the canonical mapping $SL_n(R[X])\rightarrow SL_n(R_M[X])$.
Also, by induction, we may assume $SL_n(R)=E_n(R)$ for $n\geq 3$.
Now consider the following analogue of Quillen's
theorem for elementary matrices;
\begin{quote}
Suppose $n\geq 3$ and $A\in SL_n(R[X])$. Then $A$ is realizable
over $R[X]$ if and only if $A_M\in SL_n(R_M[X])$ is realizable
over $R_M[X]$ for every $M\in$ Max($R$).
\end{quote}
While a non-constructive proof of this assertion is given in
\cite{suslin} and a more general functorial treatment of this
{\em Quillen Induction Process} can be found in \cite{knus}, we will attempt
to give a constructive proof for it here.
Since the necessity of the condition is clear, we have to prove
the following;
\begin{thm}\label{thm;glueing}
(Quillen Induction Algorithm)
For any given $A\in SL_n(R[X])$, if $A_M\in E_n(R_M[X])$
for every $M\in {\rm Max}(R)$, then  $A\in E_n(R[X])$.
\end{thm}
\begin{remark}
In view of this theorem, for any given $A\in SL_n(R[X])$,
now it's enough to have a realization
algorithm for each $A_M$ over $R_M[X]$.
\end{remark}
\pfbegin
Let ${\bf a_1}=(0,\ldots ,0)\in k^{m-1}$, and
$M_1=\{ g\in k[x_1,\ldots ,x_{m-1}]\mid g({\bf a_1})=0\}$
be the corresponding maximal ideal. Then by the condition of the
theorem, $A_{M_1}$ is realizable over $R_{M_1}[X]$. Hence, we can
write
\begin{eqnarray}
A_{M_1}=\prod_jE_{s_jt_j}\left(\frac{c_j}{d_j}\right)
\end{eqnarray}
where $c_j,d_j\in R, d_j\not\in M_1$.
Letting $r_1=\prod_jd_j\notin M_1$, we can rewrite this as
\begin{eqnarray}
A_{M_1}=\prod_jE_{s_jt_j}\left(\frac{c_j\prod_{k\neq j}d_k}{r_1}\right)
\in E_n(R_{r_1})\subset E_n(R_{M_1}).
\end{eqnarray}
Denote an algebraic closure of $k$ by $\bar{k}$.
Inductively, let ${\bf a_j} \in {\bar{k}}^{m-1}$ be a common zero of
$r_1,\ldots ,r_{j-1}$
and $M_j=\{ g\in k[x_1,\ldots ,x_{m-1}]\mid g({\bf a_j})=0\}$
be the corresponding maximal ideal of $R$ for
each $j\geq 2$. Define $r_j\notin M_j$ in the same way
as in the above so that
\begin{eqnarray}
A_{M_j}\in E_n(R_{r_j}[X]).
\end{eqnarray}
Since ${\bf a_j}$ is a common zero of
$r_1,\ldots ,r_{j-1}$ in this construction, we immediately see
$r_1,\ldots ,r_{j-1}\in M_j=\{g\in R\mid g({\bf a_j})=0\}$.
But noting $r_j\notin M_j$,
we conclude that $r_j\notin r_1R+\cdots +r_{j-1}R$.
Now, since the Noetherian condition on $R$ guarantees that we will
get to some $L$ after a finite number of steps such that
$r_1R+\cdots +r_LR=R$,
we can use the usual {\em Ideal Membership Algorithm} to determine
when $1_R$ is in the ideal $r_1R+\cdots +r_LR$.

Let $l$ be a {\em large} natural number (It will soon be clear what {\em
large} means). Then since $r_1^lR+\cdots +r_L^lR=R$, we can
use the {\em effective Nullstellensatz} to find
$g_1,\ldots ,g_L\in R$ such that $r_1^lg_1+\cdots +r_L^lg_L=1$.
Now, we express $A(X)\in SL_n(R[X])$ in the following way:
\begin{eqnarray}
A(X) & = & A(X-Xr_1^lg_1) \cdot [ A^{-1}(X-Xr_1^lg_1) A(X) ]\nonumber\\
    & = & A(X-Xr_1^lg_1-Xr_2^lg_2) \cdot
[ A^{-1}(X-Xr_1^lg_1-Xr_2^lg_2) A(X-Xr_1^lg_1) ] \nonumber\\
    &   & \cdot [ A^{-1}(X-Xr_1^lg_1) A(X) ]\nonumber\\
    & = &  \cdots \nonumber\\
    & = & A(X-\sum_{i=1}^LXr_i^lg_i) \cdot [ A^{-1}(X-\sum_{i=1}^LXr_i^lg_i)
A(X-\sum_{i=1}^{L-1}Xr_i^lg_i) ] \cdots \nonumber \\
    &   & \cdots [ A^{-1}(X-Xr_1^lg_1) A(X) ].
\end{eqnarray}
Note here that the first matrix $A(X-\sum_{i=1}^LXr_i^lg_i)=A(0)$
on the right hand side is
in $SL_n(R)=E_n(R)$ by the induction hypothesis. What will be shown
now is that for a sufficiently large $l$, each expression in the
brackets in the above equation for $A$ is actually in $E_n(R[X])$,
so that $A$ itself is in $E_n(R[X])$.
To this end, by letting $A_{M_i}=A_i$ and identifying
$A\in SL_n(R[X])$ with $A_i\in SL_n(R_{M_i}[X])$, note that each
expression in the brackets is in the following form:
\begin{eqnarray}
A_i^{-1}(cX)A_i((c+r_i^lg)X).
\end{eqnarray}
\newline {\bf *Claim:} For any $c,g\in R$, we can find a sufficiently large
$l$ such that $A_i^{-1}(cX)A_i((c+r_i^lg)X)\in E_n(R[X])$ for all
$i=1,\ldots ,L$.
\newline Let
\begin{eqnarray}
D_i(X,Y,Z)=A_i^{-1}(Y\cdot X)A_i((Y+Z)\cdot X)\in E_n(R_{r_i}[X,Y,Z])
\end{eqnarray}
and write $D_i$ in the form
\begin{eqnarray}
D_i=\prod_{j=1}^hE_{s_jt_j}(b_j+Zf_j)
\end{eqnarray}
where $b_j\in R_{r_i}[X,Y]$ and $f_j\in R_{r_i}[X,Y,Z]$.
{}From now on, the elementary matrix $E_{s_jt_j}(a)$ will be simply
denoted as $E^j(a)$ for notational convenience.
Now define $C_p$ by
\begin{eqnarray}
C_p=\prod_{j=1}^pE^j(b_j)\in E_n(R_{r_i}[X,Y]).
\end{eqnarray}
Then the $C_p$'s satisfy the following recursive relations;
\begin{eqnarray}
E^1(b_1) & = & C_1 \nonumber \\
      E^p(b_p) & = & C_{p-1}^{-1}C_p\quad (2\leq p\leq h)\nonumber \\
      C_h      & = & I.
\end{eqnarray}
Hence, using $E_{ij}(a+b)=E_{ij}(a)E_{ij}(b)$,
\begin{eqnarray}
D_i & = & \prod_{j=1}^hE^j(b_j+Zf_j) \nonumber \\
          & = & \prod_{j=1}^hE^j(b_j)E^j(Zf_j) \nonumber \\
          & = & [E^1(b_1)E^1(Zf_1)][E^2(b_2)E^2(Zf_2)]
\ \cdots\  [E^h(b_h)E^h(Zf_h)] \nonumber \\
          & = & [C_1E^1(Zf_1)][C_1^{-1}C_2E^2(Zf_2)]\ \cdots\
[C_{h-1}^{-1}C_hE^h(Zf_h)] \nonumber \\
          & = & \prod_{j=1}^hC_jE^j(Zf_j)C_j^{-1}.
\end{eqnarray}
Now in the same way as in the proof of {\bf Corollary 1} and
{\bf Corollary 2} of section 2, we can write
$C_jE^j(Zf_j)C_j^{-1}$ as a product of Cohn-type matrices, i.e.
for any given $j\in \{1,\ldots ,h\}$, let
${\bf v}=\left(\!\begin{array}{c} v_1 \\ \vdots \\ v_n \end{array}\!\right)$
be the $s_j$-th column vector of $C_j$. Then
\begin{eqnarray}
C_jE_{s_jt_j}(Zf_j)C_j^{-1}=\prod_{1\leq \g < \de \leq n}
[I+{\bf v}\cdot Zf_j\cdot a_{\g \de}(v_{\g}{\bf e_{\de}}-v_{\de}{\bf e_{\g}})]
\end{eqnarray}
for some $a_{\g \de}\in R_{r_i}[X,Y]$.
Also we can find a natural number $l$ such that
\begin{eqnarray}
v_{\g}=\frac{v_{\g}'}{r_i^l},\quad a_{\g \de}=\frac{a_{\g \de}'}{r_i^l},\quad
f_j=\frac{f_j'}{r_i^l}
\end{eqnarray}
for some $v_{\g}',a_{\g \de}'\in R[X,Y],\ f_j'\in
R[X,Y,Z]$. Now, replacing $Z$ by $r_i^{4l}g$, we see that
all the Cohn-type matrices in the above expression for $C_jE^j(Zf_j)C_j^{-1}$
have denominator-free entries. Therefore,
\begin{eqnarray}
C_jE^j(r_i^{4l}gf_j)C_j^{-1}\in E_n(R[X,Y]).
\end{eqnarray}
Since this is true for each $j$, we conclude that for a sufficiently
large $l$,
\begin{eqnarray}
D_i(X,Y,r_i^{l}g)=\prod_{j=1}^hC_jE^j(r_i^{l}gf_j)
C_j^{-1}\in  E_n(R[X,Y]).
\end{eqnarray}
Now, letting $Y=c$ proves the claim.
\pfend
\section{Reduction to $SL_3(k[x_1,\ldots ,x_m])$}
Let $A\in SL_n(k[x_1,\ldots ,x_m])$ with $n\geq 3$, and ${\bf v}$
be its last column vector. Then ${\bf v}$ is unimodular. (Recall that
the cofactor expansion along the last column gives a required relation.)
Now, if we can reduce ${\bf v}$ to ${\bf e_n}=(0,0,\ldots ,0,1)^t$
by applying elementary operations, i.e. if we can
find $B\in E_n(k[x_1,\ldots ,x_m])$ such that $B{\bf v}={\bf e_n}$,
then
\begin{eqnarray}
BA=\left(\! \begin{array}{cccc} & & & 0 \\
 & \tilde{A} & & \vdots \\
 & & & 0 \\ p_1 &\ldots  & p_{n-1} & 1 \end{array}\! \right)
\end{eqnarray}
for some $\tilde{A} \in SL_{n-1}(k[x_1,\ldots ,x_m])$ and $p_i \in
k[x_1,\ldots ,x_m]$ for
$i=1,\ldots ,n-1$. Hence,
\begin{eqnarray}
BAE_{n1}(-p_1)\cdots E_{n(n-1)}(-p_{n-1})=
\left(\!\begin{array}{cc} \tilde{A} & 0 \\ 0 & 1 \end{array}\!\right).
\end{eqnarray}
Therefore our problem of expressing $A\in SL_n(k[x_1,\ldots ,x_m])$
as a product of elementary matrices is now reduced to the same problem
for $\tilde{A}\in SL_{n-1}(k[x_1,\ldots ,x_m])$. By repeating this
process, we get to the problem of expressing $A=
\left(\! \begin{array}{ccc} p & q & 0 \\ r & s & 0 \\ 0 & 0 & 1
\end{array}\! \right) \in SL_3(k[x_1,\ldots ,x_m])$ as a product of
elementary matrices, which is the subject of the next section.
In this section, we will develop an algorithm for finding
elementary operations that reduce a given unimodular column vector
${\bf v} \in (k[x_1,\ldots ,x_m])^n$ to ${\bf e_n}$.
Also, as a corollary to this {\em Elementary Column Property}, we
give an algorithmic proof of
the {\em Unimodular Column Property} which states that
for any given unimodular column vector ${\bf v} \in
(k[x_1,\ldots ,x_m])^n$, there exists a unimodular matrix $B$,
i.e. a matrix of constant determinant, over $k[x_1,\ldots ,x_m]$
such that $B{\bf v}={\bf e_n}$. Lately, {\sl A. Logar, B. Sturmfels}
in \cite{logar:sturmfels} and {\sl N. Fitchas, A. Galligo} in
\cite{fitchas:galligo}, \cite{fitchas} have given different
algorithmic proofs of this {\em Unimodular Column Property}, thereby
giving algorithmic proofs of the {\em Quillen-Suslin theorem}.
Therefore, our algorithm gives another constructive proof of the
{\em Quillen-Suslin theorem}. The second author has given a different
algorithmic proof of the {\em Elementary Column Property} based
on a localization and patching process in \cite{cynthia}.
\begin{definition}
For a ring $R$, ${\rm Um}_n(R)=\{ n$-dimensional unimodular
column vectors over $R\}$.
\end{definition}
\begin{remark}
Note that the groups $GL_n(k[x_1,\ldots ,x_m])$
and $E_n(k[x_1,\ldots ,x_m])$ act on the set
${\rm Um}_n(k[x_1,\ldots ,x_m])$ by matrix multiplication.
\end{remark}
\begin{thm}\label{thm;reduction}
(Elementary Column Property)
For $n\geq 3$, the group $E_n(k[x_1,\ldots ,x_m])$ acts transitively
on the set Um$_n(k[x_1,\ldots ,x_m])$.
\end{thm}
\begin{remark}
According to this theorem, if ${\bf v,v'}$ are $n$-dimensional
unimodular column vectors over $k[x_1,\ldots ,x_m]$, then we can find $B\in
E_n(k[x_1,\ldots ,x_m])$
such that $B{\bf v}={\bf v'}$. Letting ${\bf v'}={\bf e_n}$
gives a desired algorithm.
\end{remark}
\begin{cor}
(Unimodular Column Property)
For $n\geq 2$, the group $GL_n(k[x_1,\ldots ,x_m])$ acts transitively
on the set Um$_n(k[x_1,\ldots ,x_m])$.
\end{cor}
\pfbegin
For $n\geq 3$, the {\em Elementary Column Property} cleary implies
the {\em Unimodular Column Property} since a product of elementary
matrices is always unimodular, i.e. has a constant determinant.

If $n=2$, for any ${\bf v}=(v_1,v_2)^t\in
\mbox{Um}_2(k[x_1,\ldots ,x_m])$,
find $g_1,g_2\in k[x_1,\ldots ,x_m]$ such that $v_1g_1+v_2g_2=1$.
Then the unimodular matrix
$U_{\bf v}= \left(\!\begin{array}{cc} v_2 & -v_1 \\ g_1 & g_2
\end{array}\!\right)$ satisfies
$U_{\bf v}\cdot {\bf v}={\bf e_2}$. Therefore we see that, for any
${\bf v}, {\bf w}\in \mbox{Um}_2(k[x_1,\ldots ,x_m])$,
$U_{\bf w}^{-1}U_{\bf v}\cdot {\bf v}={\bf w}$ where
$U_{\bf w}^{-1}U_{\bf v}\in GL_2(k[x_1,\ldots ,x_m])$.
\pfend
Let $R=k[x_1,\ldots ,x_{m-1}]$ and $X=x_m$. Then $k[x_1,\ldots ,x_m]=R[X]$.
By identifying $A\in SL_2(R[X])$ with
$\left(\! \begin{array}{cc} A & 0 \\ 0 & I_{n-2}
\end{array}\!\right)\in SL_n(R[X])$,
we can regard $SL_2(R[X])$ as a subgroup of
$SL_n(R[X])$. Now consider the following theorem.
\begin{thm}\label{thm;link}
Suppose ${\bf v}(X)=
\left(\!\begin{array}{c} v_1(X) \\ \vdots \\ v_n(X) \end{array}\!\right)
\in {\rm Um}_n(R[X])$, and $v_1(X)$ is monic in $X$.
Then there exists $B_1\in SL_2(R[X])$ and $B_2\in E_n(R[X])$ such that
$B_1B_2\cdot {\bf v}(X)={\bf v}(0)$.
\end{thm}
\pfbegin
Later
\pfend
We will use this theorem to prove the {\bf Theorem~\ref{thm;reduction}}, now.

\medskip

\noindent {\bf Proof of Theorem~\ref{thm;reduction}:}
Since the {\em Euclidean division algorithm} for $k[x_1]$ proves the
theorem for $m=1$ case,
by induction, we may assume the statement of the theorem for $R=
k[x_1,\ldots ,x_{m-1}]$. Let $X=x_m$ and
${\bf v}=\left(\!\begin{array}{c} v_1\\ \vdots\\ v_n\end{array}\!\right)
\in {\rm Um}_n(R[X])$. We may also assume that $v_1$ is monic
by applying a change of variables
(as in the well-known proof of the {\em Noether Normalization Lemma}).
Now by the above {\bf Theorem~\ref{thm;link}}, we can find
$B_1\in SL_2(R[X])$ and $B_2\in E_n(R[X])$ such that
\begin{eqnarray}
B_1B_2\cdot {\bf v}(X)={\bf v}(0)\in R.
\end{eqnarray}
And then by the inductive hypothesis,
we can find $B'\in E_n(R)$ such that
\begin{eqnarray}
B'\cdot {\bf v}(0)={\bf e_n}.
\end{eqnarray}
Therefore, we get
\begin{eqnarray}
{\bf v}=B_2^{-1}B_1^{-1}B'^{-1}{\bf e_n}.
\end{eqnarray}
By the normality of $E_n(R[X])$ in $SL_n(R[X])$ ({\bf Corollary 2}),
we can write $B_1^{-1}B'^{-1}=B''B_1^{-1}$ for some $B''\in E_n(R[X])$.
Since
\begin{eqnarray}
B_1^{-1}=\left(\! \begin{array}{ccccc}
p & q& 0 &\ldots  & 0 \\ r & s & 0 & \ldots &0\\
0 & 0 &&& \\\vdots & \vdots &  & I_{n-2} & \\ 0 & 0 &&&
\end{array}\!\right)
\end{eqnarray}
for some $p,q,r,s \in R[X]$, we have
\begin{eqnarray}
{\bf v} & = & B_2^{-1}B_1^{-1}B'^{-1}{\bf e_n}\nonumber\\
	& = & (B_2^{-1}B'')B_1^{-1}{\bf e_n}\nonumber\\
	& = & (B_2^{-1}B'')\left(\! \begin{array}{ccccc}
p & q& 0 &\ldots  & 0 \\ r & s & 0 & \ldots &0\\
0 & 0 &&& \\\vdots & \vdots &  & I_{n-2} & \\ 0 & 0 &&&
\end{array}\!\right)\left(\! \begin{array}{c} 0 \\ 0 \\ \vdots
\\ 0 \\ 1 \end{array} \!\right)\nonumber\\
	& = & (B_2^{-1}B''){\bf e_n}
\end{eqnarray}
where $B_2^{-1}B''\in E_n(R[X])$.
Since we have this relationship for any ${\bf v}\in {\rm Um}_n(R[X])$,
we get the desired transitivity.
\pfend
Now, we need one lemma to construct an algorithm for
the {\bf Theorem~\ref{thm;link}}.
\begin{lemma}\label{lem;link}
Let $f_1,f_2,b,d\in R[X]$ and $r$ be the resultant of $f_1$ and $f_2$.
Then there exists $B\in SL_2(R[X])$ such that
\begin{eqnarray}
B\left(\! \begin{array}{c} f_1(b) \\f_2(b) \end{array}\!\right)
=\left( \!\begin{array}{c} f_1(b+rd) \\f_2(b+rd) \end{array}\!\right).
\end{eqnarray}
\end{lemma}
\pfbegin
By the property of the resultant of two polynomials, we can find
$g_1,g_2\in R[X]$ such that $f_1g_1+f_2g_2=r$. Also let
$s_1,s_2,t_1,t_2\in R[X,Y,Z]$ be the polynomials defined by
\begin{eqnarray}
f_1(X+YZ) & = & f_1(X)+Ys_1(X,Y,Z)\nonumber\\
f_2(X+YZ) & = & f_2(X)+ Ys_2(X,Y,Z)\nonumber\\
g_1(X+YZ) & = & g_1(X)+Yt_1(X,Y,Z) \nonumber\\
g_2(X+YZ) & = & g_2(X)+ Yt_2(X,Y,Z).
\end{eqnarray}
Now, let
\begin{eqnarray}
B_{11} & = & 1+s_1(b,r,d)\cdot g_1(b)+t_2(b,r,d)\cdot f_2(b) \nonumber\\
B_{12} & = & s_1(b,r,d)\cdot g_2(b)-t_2(b,r,d)\cdot f_1(b) \nonumber \\
B_{21} & = & s_2(b,r,d)\cdot g_1(b)-t_1(b,r,d)\cdot f_2(b) \nonumber \\
B_{22} & = & 1+s_2(b,r,d)\cdot g_2(b)+t_1(b,r,d)\cdot f_1(b).
\end{eqnarray}
Then one checks easily that
$B=\left(\! \begin{array}{cc} B_{11} & B_{12} \\ B_{21} & B_{22}
\end{array}\!\right)$ satisfies the desired property and that
$B\in SL_2(R[X])$.
\pfend
{\bf Proof of Theorem~\ref{thm;link}:}
Let ${\bf a_1}=(0,\ldots ,0)\in k^{m-1}$. Define
$M_1=\{ g\in k[x_1,\ldots ,x_{m-1}]\mid g({\bf a_1})=0\} $
and $k_1=R/M_1$ as the corresponding maximal ideal and residue
field, respectively.
Since ${\bf v}\in (R[X])^n$ is a unimodular column vector, its image
${\bf \bar{v}}$ in $(k_1[X])^n=((R/M_1)[X])^n$ is also unimodular.
Since $k_1[X]$ is a principal ideal ring, the minimal Gr\"{o}bner basis of
its ideal $<\bar{v}_2,\ldots ,\bar{v}_n>$ consists of
a single element, $G_1$.
Then $\bar{v}_1$ and $G_1$ generate the unit ideal in $k_1[X]$ since
$\bar{v}_1,\bar{v}_2,\ldots ,\bar{v}_n$ generate the unit ideal.
Using the Euclidean division algorithm for $k_1[X]$, we can find $E_1\in
E_{n-1}(k_1[X])$ such that
\begin{eqnarray}
E_1\left( \!\begin{array}{c} \bar{v}_2  \\ \vdots \\ \bar{v}_n
\end{array} \right)=
\left( \begin{array}{c} G_1 \\ 0 \\ \vdots \\ 0
\end{array}\! \right).
\end{eqnarray}
By identifying $k_1$ with a subring of $R$, we may regard
$E_1$ to be an element of $E_n(R[X])$
and $G_1$ to be an element of $R[X]$.
Then,
\begin{eqnarray}
\left( \!\begin{array}{cc} 1 & 0 \\ 0 & E_1 \end{array} \right)
{\bf v}=\left( \begin{array}{c} v_1
\\ G_1+q_{12} \\ q_{13}
\\ \vdots \\ q_{1n}
\end{array} \!\right)
\end{eqnarray}
for some $q_{12},\ldots, q_{1n}\in M_1[X]$.
Now, define $r_1\in R$ by
\begin{eqnarray}
r_1 & = & {\rm Res}(v_1, G_1+q_{12}) \nonumber\\
	& = & {\rm the\  resultant\  of}\
v_1\ {\rm and}\ G_1+q_{12}
\end{eqnarray}
and find $f_1,h_1\in R[X]$ such that
\begin{eqnarray}
f_1\cdot v_1+h_1\cdot (G_1+q_{12})=r_1.
\end{eqnarray}
Since $v_1$ is monic, and $\bar{v}_1$ and $G_1\in k_1[X]$ generate the unit
ideal, we have
\begin{eqnarray}
\bar{r}_1 & = & \overline{{\rm Res}(v_1, G_1+q_{12})}\nonumber\\
	& = & {\rm Res}(\bar{v}_1, G_1)\nonumber\\
	& \neq & 0.
\end{eqnarray}
Therefore, $r_1\notin M_1$. Denote an algebraic closure of $k$ by
$\bar{k}$.
Inductively, let ${\bf a_j}\in {\bar{k}}^{m-1}$ be a common zero of
$r_1,\ldots ,r_{j-1}$
and $M_j$ be the corresponding maximal ideal of $R$ for
each $j\geq 2$. Define $r_j\notin M_j$ in the same way
as in the above. Define also, $E_j\in E_{n-1}(k_j[X]), G_j\in k_j[X],
f_j,h_j\in R[X]$,
and $q_{j2}, \ldots , q_{jn}\in M_j[X]$ in an analogous way.
Since we let ${\bf a_j}$ be a common zero of
$r_1,\ldots ,r_{j-1}$ in this construction, we see
$r_1,\ldots ,r_{j-1}\in M_j=\{g\in R\mid g({\bf a_j})=0\}$.
But noting $r_j\notin M_j$,
we conclude that $r_j\notin r_1R+\cdots +r_{j-1}R$.
Now, since $R$ is Noetherian, after a
finite number of steps, we will get to some $L$ such that
$r_1R+\cdots +r_LR=R$.
We can use the {\em effective Nullstellensatz} to explicitly find
those $g_i$'s in
$R$ such that $r_1g_1+\cdots +r_Lg_L= 1$.
Define, now, $b_0,b_1,\ldots ,b_L\in R[X]$ in the following way:
\begin{eqnarray}
	b_0 & = & 0 \nonumber\\
	b_1 & = & r_1g_1X \nonumber\\
	b_2 & = & r_1g_1X+r_2g_2X \nonumber\\
	   & \vdots & \nonumber\\
	b_L & = & r_1g_1X+r_2g_2X+\cdots +r_Lg_LX=X.
\end{eqnarray}
Then these $b_i$'s satisfy the recursive relations:
\begin{eqnarray}
	b_0 & = & 0 \nonumber\\
	b_i & = & b_{i-1} + r_ig_iX \quad {\rm for}\ i=1,\ldots ,L.
\end{eqnarray}
{\bf *Claim:}
For each $i\in \{ 1,\ldots ,L\}$, there exists $B_i\in SL_2(R[X])$
and $B_i'\in E_n(R[X])$ such that
${\bf v}(b_i)=B_iB_i'{\bf v}(b_{i-1})$.

If this claim is true, then using $E_n(R[X])\cdot SL_2(R[X])\subseteq
SL_2(R[X])\cdot E_n(R[X])$ (Normality of $E_n(R[X])$; {\bf Corollary 2}),
we inductively get
\begin{eqnarray}
{\bf v}(X) & = & {\bf v}(b_L)\nonumber\\
	& = & B_LB_L'{\bf v}(b_{L-1})\nonumber\\
	& \vdots & \nonumber\\
	& = & BB'{\bf v}(b_0)\nonumber\\
	& = & BB'{\bf v}(0)
\end{eqnarray}
for some $B\in SL_2(R[X])$
and $B'\in E_n(R[X])$. Therefore it's enough to prove the above
claim.
For this purpose, let $\tilde{G_i} =G_i+q_{i2}$. Then
\begin{eqnarray}
\left(\! \begin{array}{cc} 1 & 0 \\ 0 & E(X) \end{array} \!\right)
{\bf v}(X)=\left(\! \begin{array}{c} v_1(X) \\ \tilde{G_i}(X) \\
q_{i3}(X) \\ \vdots \\ q_{in}(X) \end{array} \!\right).
\end{eqnarray}
For $3\leq l \leq n$,
we have
\begin{eqnarray}
q_{il}(b_i)-q_{il}(b_{i-1}) & \in & (b_i - b_{i-1})\cdot R[X] \nonumber\\
		& = &  r_ig_iX \cdot R[X].
\end{eqnarray}
Since $r_i\in R$ doesn't depend on $X$, we have
\begin{eqnarray}
r_i & = & f_i(X)v_1(X)+h_i(X)\tilde{G_i}(X)\nonumber\\
	& = & f_i(b_{i-1})v_1(b_{i-1})+h_i(b_{i-1})\tilde{G_i}
(b_{i-1})\nonumber\\
	& = & {\rm a\ linear\ combination\ of}\ v_1(b_{i-1})\
{\rm and}\ \tilde{G_i}(b_{i-1})\ {\rm over}\ R[X].
\end{eqnarray}
Therefore, we see that for $3\leq l \leq n$,
\begin{eqnarray}
q_{il}(b_i)=q_{il}(b_{i-1})+{\rm a\ linear\ combination\ of}\ v_1(b_{i-1})\
{\rm and}\ \tilde{G_i}(b_{i-1})\ {\rm over} R[X].\nonumber
\end{eqnarray}
Hence we can find $C\in E_n(R[X])$ such that
\begin{eqnarray}
C \left(\! \begin{array}{cc} 1 & 0 \\ 0 & E(b_{i-1}) \end{array} \!\right)
{\bf v}(b_{i-1}) & = &
C \left(\! \begin{array}{c} v_1(b_{i-1}) \\ \tilde{G}_i(b_{i-1}) \\
q_{i3}(b_{i-1}) \\ \vdots \\ q_{in}(b_{i-1}) \end{array} \!\right)
\nonumber\\
	& = & \left(\! \begin{array}{c} v_1(b_{i-1})\\ \tilde{G}_i(b_{i-1}) \\
q_{i3}(b_{i}) \\ \vdots \\ q_{in}(b_{i}) \end{array} \!\right).
\end{eqnarray}
Now, by the {\bf Lemma~\ref{lem;link}}, we can find
$\tilde{B} \in SL_2(R[X])$ such that
\begin{eqnarray}
\tilde{B} \left(\! \begin{array}{c} v_1(b_{i-1}) \\ \tilde{G_i}(b_{i-1})
\end{array} \!\right)=\left(\! \begin{array}{c} v_1(b_i) \\ \tilde{G_i}(b_i)
\end{array} \!\right).
\end{eqnarray}
Finally, define $B\in SL_n(R[X])$ as follows:
\begin{eqnarray}
B=\left(\! \begin{array}{cc} 1 & 0 \\ 0 & E(b_i)^{-1} \end{array} \!\right)
\left(\! \begin{array}{cc} \tilde{B} & 0 \\ 0 & I_{n-2} \end{array} \!\right)
\cdot C\cdot \left(\!\begin{array}{cc} 1 &0 \\ 0& E(b_i)\end{array}\!\right).
\end{eqnarray}
Then this $B$ satisfies
\begin{eqnarray}
B{\bf v}(b_{i-1})={\bf v}(b_i),
\end{eqnarray}
and by using the normality of $E_n(R[X])$ again, we see that
\begin{eqnarray}
B\in SL_2(R[X])E_n(R[X])
\end{eqnarray}
and this proves the claim.
\pfend
\section{Realization Algorithm for $SL_3(R[X])$}
Now, we want to find a realization
algorithm for the matrices of the special type in $SL_3(k[x_1,\ldots ,x_m])$,
i.e.
matrices of the form $\left(\!\begin{array}{ccc} p & q & 0 \\ r & s & 0 \\ 0 &
0 & 1
\end{array}\!\right) \in SL_3(k[x_1,\ldots ,x_m])$. Again, by applying
a change of variables, we may assume
that $p\in k[x_1,\ldots ,x_m]$ is a monic polynomial in the
last variable $x_m$. In view of the
{\em Quillen Induction Algorithm} developed in the section 3, we see that
it's enough to develop a realization algorithm for the matrices of the
form $\left(\!\begin{array}{ccc} p & q & 0 \\ r & s & 0 \\ 0 & 0 & 1
\end{array}\! \right) \in SL_3(R[X])$, where $R$ is now a commutative local
ring and $p\in R[X]$ is a monic polynomial. A realization algorithm for this
case
was obtained by M.P. Murthy, and we present in the below a slightly
modified version of the {\bf Lemma 3.6}
in \cite{murthy} {\sl Suslin's Work on Linear Groups over Polynomial Rings
and Serre Problem} by S.K. Gupta and M.P. Murthy.
\begin{lemma}\label{lemma:split}
Let $L$ be a commutative ring, and $a,a',b\in L$. Then, the followings are
true.
\begin{enumerate}
\item $(a,b)$ and $(a',b)$ are unimodular over $L$ if and only if
$(aa',b)$ is unimodular over $L$.
\item For any $c,d\in L$ such that $aa'd-bc=1$,
there exist $c_1,c_2,d_1,d_2\in L$ such that $ad_1-bc_1=1,\ a'd_2-bc_2=1$, and
$$\left(\!\begin{array}{ccc} aa' & b & 0 \\ c & d & 0 \\ 0 & 0 & 1
\end{array}\!\right)\equiv \left(\!\begin{array}{ccc}a&b&0\\ c_1&d_1&0\\ 0&0&1
\end{array}\!\right)\cdot \left(\!\begin{array}{ccc}a'&b&0\\ c_2&d_2&0\\ 0&0&1
\end{array}\!\right) \pmod {E_3(L)}.$$
\end{enumerate}
\end{lemma}
\pfbegin
(1) If $(aa',b)$ is unimodular over $L$, there exist $h_1,h_2\in L$ such that
$h_1\cdot (aa')+h_2\cdot b=1$. Now $(h_1a')\cdot a+h_2\cdot b=1$ implies
$(a,b)$ is unimodular, and $(h_1a)\cdot a'+h_2\cdot b=1$ implies $(a',b)$
is unimodular.

Suppose, now, that $(a,b)$ and $(a',b)$ are unimodular over $L$. Then, we can
find $h_1,h_2,h_1',h_2'\in L$ such that
$h_1a+h_2b=1,\ h_1'a'+h_2'b=1$.
Now, let $g_1=h_1h_1',\ g_2=h_2'+a'h_2h_1'$, and consider
\begin{eqnarray}
g_1aa'+g_2b & = & h_1h_1'aa'+(h_2'+a'h_2h_1')b \nonumber\\
            & = & h_1'a'(h_1a+h_2b)+h_2'b\nonumber\\
            & = & h_1'a'+h_2'b\nonumber\\
            & = & 1.
\end{eqnarray}
So we have a desired unimodular relation.

\medskip
\noindent (2) If $c,d\in L$ satisfy $aa'd-bc=1$, then $(aa',b)$ is
unimodular, which in turn implies that $(a,b)$ and $(a',b)$ are unimodular.
Therefore, we can find
$c_1,d_1,d_1,d_2\in L$ such that  $ad_1-bc_1=1$ and $a'd_2-bc_2=1$.
For example, we can let
\begin{eqnarray}
c_1=c_2=c,\quad d_1=a'd, \quad d_2=ad.
\end{eqnarray}
Now, consider
\begin{eqnarray}
\left(\!\begin{array}{ccc} aa' & b & 0 \\ c & d & 0 \\ 0 & 0 & 1
\end{array}\!\right) & = & E_{21}(cd_1d_2-d(c_2+a'c_1d_2))
\left(\!\begin{array}{ccc}aa'&b& 0\\ c_2+a'c_1d_2&d_1d_2 &0 \\ 0&0&1
\end{array}\!\right)\nonumber\\
& = & E_{21}(cd_1d_2-d(c_2+a'c_1d_2))E_{23}(d_2-1)E_{32}(1)E_{23}(-1)
\nonumber\\
& & \left(\!\begin{array}{ccc} a & b & 0\\ c_1 & d_1 & 0\\ 0 & 0 & 1
\end{array}\!\right)E_{23}(1)E_{32}(-1)E_{23}(1)
\left(\!\begin{array}{ccc}a'&b&0\\ c_2&d_2&0 \\ 0&0&1\end{array}
\!\right)\nonumber\\
& & E_{23}(-1)E_{32}(1)E_{23}(a-1)E_{31}(-a'c_1)E_{32}(-d_1).
\end{eqnarray}
This explicit expression tells us that
\begin{eqnarray}
\left(\!\begin{array}{ccc} aa' & b & 0 \\ c & d & 0 \\ 0 & 0 & 1
\end{array}\!\right) & \equiv & \left(\!\begin{array}{ccc} a&b&0\\ c_1&d_1&0
\\ 0&0&1\end{array}\!\right) \cdot \left(\!\begin{array}{ccc}a'&b&0
\\ c_2&d_2&0\\ 0&0&1\end{array}\!\right)\pmod {E_3(L)}.
\end{eqnarray}
\pfend
\begin{thm}
Suppose $(R,M)$ is a commutative local ring, and
$A=\left(\!\begin{array}{ccc} p & q & 0 \\ r & s & 0 \\ 0 & 0 & 1
\end{array}\!\right) \in SL_3(R[X])$ where $p$ is monic.
Then $A$ is realizable over $R[X]$.
\end{thm}
\pfbegin
By induction on $\deg (p)$. If $\deg (p)=0$, then $p=0\ \mbox{or}\ 1$, and
$A$ is clearly
realizable. Now, suppose $\deg(p)=d>0$ and $\deg(q)=l$. Since $p\in R[X]$
is monic, we can find $f,g\in R[X]$ such that
\begin{eqnarray}
q & = & fp+g,\quad \deg(g)<d.
\end{eqnarray}
Then,
\begin{eqnarray}
AE_{12}(-f) & = & \left(\!\begin{array}{ccc}p&q-fp&0\\ r&s-fr&0\\ 0&0 & 1
\end{array}\!\right)=\left(\!\begin{array}{ccc}p&g&0\\ r&s-fr&0\\ 0&0 & 1
\end{array}\!\right).
\end{eqnarray}
Hence we may assume $\deg(q)<d$. Now, we note that either $p(0)$ or
$q(0)$ is a unit in $R$, otherwise, we would have $p(0)s(0)-q(0)r(0)\in M$
that contradicts to $ps-qr=p(0)s(0)-q(0)r(0)=1$. Let's consider these
two cases, separately.

\medskip
\noindent Case 1: When $q(0)$ is a unit.
\newline Using the invertibility of $q(0)$, we have
\begin{eqnarray}
AE_{21}(-q(0)^{-1}p(0)) & = & \left(\!\begin{array}{ccc}
p-q(0)^{-1}p(0)q & q&0\\ r-q(0)^{-1}p(0)s& s &0\\ 0&0&1
\end{array}\!\right).
\end{eqnarray}
So, we may assume $p(0)=0$. Now, write $p=Xp'$. Then, by the above
{\bf Lemma}~\ref{lemma:split}, we can find $c_1,d_1,c_2,d_2\in R[X]$
such that $Xd_1-qc_1=1,\ p'd_2-qc_2=1$ and
\begin{eqnarray}
\left(\!\begin{array}{ccc} p & q & 0 \\ r & s & 0 \\ 0 & 0 & 1
\end{array}\!\right) & \equiv & \left(\!\begin{array}{ccc}X&q&0\\ c_1&d_1&0
\\ 0&0&1\end{array}\!\right) \cdot \left(\!\begin{array}{ccc}p'&q&0
\\ c_2&d_2&0\\ 0&0&1\end{array}\!\right)\pmod {E_3(R[X])}
\end{eqnarray}
Since $\deg(p')<d$, the second matrix on the right hand side is
realizable by the induction hypothesis. As for the first one, we may
assume that $q$ is a unit of $R$ since we can assume $\deg(q)<\deg(X)=1$
and $q(0)$ is a unit. And then invertibility of $q$ leads easily to
an explicit factorization of $\left(\!\begin{array}{ccc}X&q&0\\ c_1&d_1&0
\\ 0&0&1\end{array}\!\right)$ into elementary matrices.

\medskip
\noindent Case 2: When $q(0)$ is not a unit.
\newline First we claim the following; there exist $p',q'\in R[X]$
such that $\deg(p')<l,\deg(q')<d$ and $p'p-q'q=1$. To prove this
claim, we let $r\in R$ be the resultant of $p$ and $q$. Then, there
exist $f,g\in R[X]$ with $\deg(f)<l,\deg(g)<d$ such that
$fp+gq=r$. Since $p$ is monic and $p,q\in R[X]$ generate the unit ideal,
we see that $r\notin M$, i.e. $r\in A^{*}$. Now, letting
$p'=f/r,\ q'=-g/r$ shows the claim. Also note that the two relations,
$p'(0)p(0)-q'(0)q(0)=1$ and $q(0)\in M$, imply $p'(0)\notin M$. This means
$q(0)+p'(0)$ is a unit. Now, consider the following.
\begin{eqnarray}
\left(\!\begin{array}{ccc} p & q & 0 \\ r & s & 0 \\ 0 & 0 & 1
\end{array}\!\right) & = & E_{21}(rp'-sq')
\left(\!\begin{array}{ccc} p & q & 0 \\ q' & p' & 0 \\ 0 & 0 & 1
\end{array}\!\right)\nonumber\\
& = & E_{21}(rp'-sq')E_{12}(-1)
\left(\!\begin{array}{ccc} p+q'& q+p'&0\\ q' & p' & 0 \\ 0 & 0 & 1
\end{array}\!\right).
\end{eqnarray}
Noting that the last matrix on the right hand side is realizable
by the Case~1 since $q(0)+p'(0)$ is a unit and $\deg(p+q')=d$, we see
that $\left(\!\begin{array}{ccc} p & q & 0 \\ r & s & 0 \\ 0 & 0 & 1
\end{array}\!\right)$ is also realizable.
\pfend
\section{Eliminating Redundancies}
When applied to a specific matrix, the algorithm presented in this
paper will produce a factorization into elementary matrices, but
this factorization may contain more factors than is necessary.
The {\em Steinberg relations} \cite{milnor}
from algebraic $K$--theory provide a method for
improving  a given  factorization by eliminating some of the
unnecessary factors.
The  {\em Steinberg relations} that elementary matrices satisfy are
\begin{enumerate}
\item
$E_{ij}(0) = I$
\item
$E_{ij}(a)E_{ij}(b) = E_{ij}(a+b)$
\item
For $i\neq l$, $[E_{ij}(a),E_{jl}(b)]=E_{ij}(a)E_{jl}(b)E_{ij}(-a)E_{jl}(-b)
= E_{il}(ab)$
\item
For $j\neq l$, $[E_{ij}(a),E_{li}(b)]= E_{ij}(a)E_{li}(b)E_{ij}(-a)E_{li}(-b)
= E_{lj}(-ab)$
\item
For $i\neq p$, $j\neq l$, $[E_{ij}(a),E_{lp}(b)]= E_{ij}(a)E_{lp}(b)E_{ij}(-a)
E_{lp}(-b)=I.$
\end{enumerate}
The first author is in the process of implementing the realization algorithm
of this paper, together with a {\em Redundancy Elimination Algorithm}
based on the above set of relations, using existing computer algebra systems.
As suggested in \cite{thk}, an algorithm of this kind has application
in {\em Signal Processing} since it gives a way of expressing a given
multidimensional filter bank as a cascade of simpler filter banks.
\section{Acknowledgement}
The authors wish to thank A. Kalker, T.Y. Lam, R. Laubenbacher, B. Sturmfels
and M. Vetterli for all the valuable support, insightful discussions
and  encouragement.

\end{document}